\documentstyle[11pt]{article}
  \newlength{\standardunitlength}
 \newtheorem{lemma}{Lemma}
\newtheorem{theorem}{Theorem} 
\newenvironment{proof}{\noindent {\sc Proof:}}{$\Box$ \vspace{2 ex}}

\begin{document}

\begin{center}
The Rogers-Ramanujan Identities, the Finite General Linear Groups, and the
Hall-Littlewood Polynomials
\end{center}

\begin{center}
By Jason Fulman
\end{center}

\begin{center}
Dartmouth College
\end{center}

\begin{center}
Jason.E.Fulman@Dartmouth.edu
\end{center}

\begin{abstract}
	The Rogers-Ramanujan identities have been studied from the
viewpoints of combinatorics, number theory, affine Lie algebras,
statistical mechanics, and quantum field theory. This note connects the
Rogers-Ramanujan identities with the finite general linear groups and the
Hall-Littlewood polynomials of symmetric function theory.
\end{abstract}

\section{Introduction}

	The Rogers-Ramanujan identities are among the most famous partition
identities in all of number theory and combinatorics. This note will be
concerned with the following generalization of the Rogers-Ramanujan
identities, due to Gordon. Let $(x)_n$ denote $(1-x)(1-x^2) \cdots
(1-x^n)$.

\begin{theorem} \label{Gordon} (Andrews \cite{Andpart}, page 111) For
$1\leq i \leq k, k \geq 2$, and $|x|<1$

\[ \sum_{n_1,\cdots,n_{k-1} \geq 0} \frac{x^{N_1^2 + \cdots + N_{k-1}^2+
N_i + \cdots+ N_{k-1}}}{(x)_{n_1}\cdots(x)_{n_{k-1}}} = \prod_{r=1 \atop r
\neq 0, \pm i (mod \ 2k+1)}^{\infty} \frac{1}{1-x^r} \]

	where $N_j = n_j + \cdots n_{k-1}$.

\end{theorem}

	Gordon's generalization of the Rogers-Ramanujan identities have
been widely studied and appear in many places in mathematics and
physics. Andrews \cite{Andpart} discusses combinatorial aspects of these
identities. Berndt \cite{Berndt} describes some number theoretic
connections. Feigen and Frenkel \cite{Feigen} interpret the Gordon
identities as a character formula for the Virasoro algebra. Andrews,
Baxter, and Forrester \cite{Baxter} relate Gordon's generalization to
exactly solved models in statistical mechanics. Berkovich, McCoy, and
Orrick \cite{Berk} and Kirillov \cite{Kirillov1}, \cite{Kirillov2} discuss
the connection with quantum field theory.

	This note studies the $i=k$ case of Gordon's generalization of the
Rogers-Ramanujan identities. Section \ref{GEN} describes the relation with
the finite general linear groups. Section \ref{SYM} describes the relation
with symmetric function theory. The results here are taken from the
Ph.D. thesis of Fulman \cite{fulthesis}.

	We use the following standard notation from the theory of
partitions. $\lambda$ is said to be a partition of $n=|\lambda|$ if
$\lambda_1 \geq \lambda_2 \geq \cdots \geq 0$ and $\sum_i \lambda_i=n$. The
$\lambda_i$ are referred to as the parts of $\lambda$. Let $m_i(\lambda)$
be the number of parts of $\lambda$ of size $i$, and define $\lambda_i'=
m_i(\lambda) + m_{i+1}(\lambda) + \cdots$. By $n(\lambda)$ is meant
$\sum_{i \geq 1} (i-1) \lambda_i$.

\section{Relation with the Finite General Linear Groups} \label{GEN}

	Recall from elementary linear algebra (for instance from Chapter 6
of Herstein \cite{Her}) that the conjugacy classes of $GL(n,q)$ are
parameterized by rational canonical form. This form corresponds to the
following combinatorial data. To each monic non-constant irreducible
polynomial $\phi$ over a field of $q$ elements, associate a partition
(perhaps the trivial partition) $\lambda_{\phi}$ of some non-negative
integer $|\lambda_{\phi}|$. Let $m_{\phi}$ denote the degree of $\phi$. The
only restrictions necessary for this data to represent a conjugacy class
are:

\begin{enumerate}

\item $|\lambda_z| = 0$
\item $\sum_{\phi} |\lambda_{\phi}| m_{\phi} = n$

\end{enumerate}

{\bf Definition} For $\alpha \in GL(n,q)$ and $\phi$ a monic, irreducible
polynomial over $F_q$, a field of $q$ elements, define
$\lambda_{\phi}(\alpha)$ to be the partition corresponding to the
polynomial $\phi$ in the rational canonical form of $\alpha$.

	The following elementary lemmas will be of use.

\begin{lemma} \label{bign} If $f(1)<\infty$ and $f$ has a Taylor
series around 0, then

\[ lim_{n \rightarrow \infty} [u^n] \frac{f(u)}{1-u} = f(1) \]

\end{lemma}

\begin{proof}
	Write the Taylor expansion $f(u) = \sum_{n=0}^{\infty} a_n
u^n$. Then observe that $[u^n] \frac{f(u)}{1-u} = \sum_{i=0}^n a_i$.
\end{proof}

	Lemma \ref{product} is proved by Stong \cite{St1} (using the fact that there are
$q^{n(n-1)}$ unipotent elements in $GL(n,q)$), but we give a simpler proof.

\begin{lemma} \label{product}

\[ 1-u = \prod_{\phi \neq z} \prod_{r=1}^{\infty}
(1-\frac{u^{m_{\phi}}}{q^{m_{\phi}r}}) \]

\end{lemma}

\begin{proof}
	First we claim that:

\[ \prod_{\phi} (1-\frac{u^{m_{\phi}}}{q^{m_{\phi}t}}) =
1-\frac{u}{q^{t-1}} \]

	Assume that $t=1$, the general case following by replacing $u$
by $\frac{u}{q^{t-1}}$. Expanding $\frac{1}{1-\frac{u^{m_{\phi}}}
{q^{m_{\phi}}}}$ as a geometric series, the coefficient of $u^d$ in
the reciprocal of the left hand side is $\frac{1}{q^d}$ times the
number of monic polynomials of degree $d$, hence 1. Comparing with the
reciprocal of the right hand side proves the claim.

	Therefore,

\[ \prod_{\phi} \prod_{r \geq 1} (1-\frac{u^{m_{\phi}}}{q^{rm_{\phi}}}) =
\prod_{r \geq 1} (1-\frac{u}{q^{r-1}}) \]

	The result follows by cancelling the terms corresponding to $\phi=z$.
\end{proof}
	
	Theorem \ref{Rogers} relates the Gordon identities with the
finite general linear groups.

\begin{theorem} \label{Rogers} Let $\phi$ be a monic, irreducible
polynomial over $F_q$ of degree $m_{\phi}$. Let $k \geq 2$ be an
integer. Then the $n \rightarrow \infty$ limit of the chance that a
uniformly chosen element $\alpha$ of $GL(n,q)$ has the largest part of the
partition $\lambda_{\phi}(\alpha)$ less than $k$ is equal to:

\[ \prod_{r=1 \atop r=0, \pm k (mod \ 2k+1)}^{\infty}
(1-\frac{1}{q^{m_{\phi}r}}) \]

\end{theorem}

\begin{proof}
	Assume for simplicity of notation that $\phi=z-1$. From the proof
it will be clear that the general case follows.	

	Stong \cite{St1}, using Kung's \cite{Kun} formula for the sizes of
the conjugacy classes of $GL(n,q)$, established the following ``cycle
index'' for the general linear groups:

\[ 1 + \sum_{n=1}^{\infty} \frac{u^n}{|GL(n,q)|} \sum_{\alpha \in GL(n,q)}
\prod_{\phi \neq z} x_{\phi,\lambda_{\phi}(\alpha)} = \prod_{\phi \neq z}
[\sum_{\lambda} x_{\phi,\lambda} \frac{u^{|\lambda|m_{\phi}}} {\prod_i
\prod_{k=1}^{m_i} (q^{m_{\phi}d_i} - q^{m_{\phi}(d_i-k)})}] \]

	where the product is over all monic $\phi \neq z$ which are
irreducible polynomials over the field of $q$ elements, and 

\[ d_i(\lambda)= 1 m_1(\lambda) + 2 m_2 (\lambda) + \cdots + (i-1)
m_{i-1}(\lambda) + (m_i(\lambda) + m_{i+1}(\lambda) + \cdots +
m_j(\lambda)) i \]

	Observe that:

\begin{eqnarray*}
\prod_i \prod_{k=1}^{m_i} (q^{m_{\phi}d_i} - q^{m_{\phi}(d_i-k)}) & = &
q^{m_{\phi} \sum_i m_i(\lambda) d_i(\lambda)} \prod_i (\frac{1}{q^{m_{\phi}}})_{m_i(\lambda)}\\
& = & q^{m_{\phi} \sum_i m_i(\lambda) [(\sum_{h<i} hm_h(\lambda)) + im_i(\lambda) +
(\sum_{i<k} im_k(\lambda))]} \prod_i (\frac{1}{q^{m_{\phi}}})_{m_i(\lambda)}\\
& = & q^{m_{\phi} \sum_i [i m_i(\lambda)^2 + 2 m_i(\lambda) \sum_{h<i} h
m_h(\lambda)]} \prod_i (\frac{1}{q^{m_{\phi}}})_{m_i(\lambda)}\\
& = & q^{m_{\phi} \sum_i (\lambda_i')^2} \prod_i (\frac{1}{q^{m_{\phi}}})_{m_i(\lambda)}
\end{eqnarray*}

	Combining this observation with Lemma \ref{product} shows that:

\[ 1 + \sum_{n=1}^{\infty} \frac{u^n}{|GL(n,q)|} \sum_{\alpha \in GL(n,q)}
\prod_{\phi \neq z} x_{\phi,\lambda_{\phi}(\alpha)} = \frac{1}{1-u}
\prod_{\phi \neq z} \prod_{r=1}^{\infty} (1-\frac{u^{m_{\phi}}}{q^{m_{\phi}r}})
[\sum_{\lambda} x_{\phi,\lambda} \frac{u^{|\lambda|m_{\phi}}} { q^{m_{\phi}
\sum_i (\lambda_i')^2} \prod_i (\frac{1}{q^{m_{\phi}}})_{m_i(\lambda)}}] \]

	Setting $x_{z-1,\lambda}=0$ if the largest part of $\lambda$ is
greater than equal to $k$, and all $x_{\phi,\lambda}=1$ otherwise shows
that the sought limiting probability is (where $[u^n]$ denotes the
coefficient of $u^n$):

\begin{eqnarray*}
& & \lim_{n \rightarrow \infty} [u^n] \frac{1}{1-u} (\prod_{r=1}^{\infty}
(1-\frac{u}{q^r})) (\sum_{\lambda: \lambda_1<k} \frac{u^{|\lambda|}} { q^{
\sum_i (\lambda_i')^2} \prod_i (\frac{1}{q})_{m_i(\lambda)}})\\
& = & (\prod_{r=1}^{\infty} (1-\frac{1}{q^r})) (\sum_{\lambda: \lambda_1<k}
\frac{1} { q^{ \sum_i (\lambda_i')^2} \prod_i
(\frac{1}{q})_{m_i(\lambda)}})\\
& = & \prod_{r=1 \atop r=0, \pm k (mod \ 2k+1)}^{\infty}
(1-\frac{1}{q^r})
\end{eqnarray*}

	The first equality uses Lemma \ref{bign} and the second equality
uses the Gordon's generalization of the Rogers-Ramanujans identities with
$n_i=m_i(\lambda)$, $i=k$, and $x=\frac{1}{q}$.
\end{proof}

{\bf Remarks}

\begin{enumerate}

\item Recall that an $\alpha \in Mat(n,q)$, all $n*n$ matrices with entries in
$F_q$, is said to be semisimple if it is diagonalizable over $\bar{F_q}$, the
algebraic closure of
$F_q$. It is elementary to show that $\alpha$ is semisimple if and only if all
$\lambda_{\phi}(\alpha)$ have largest part at most one. Combining this with
Theorem \ref{Rogers} and Stong's cycle index for $Mat(n,q)$ one can prove that
the $n \rightarrow \infty$ limiting probability that an element of
$Mat(n,q)$ is semisimple is:

\[ \prod_{r=1 \atop r=0,\pm 2 (mod \ 5)}^{\infty} (1-\frac{1}{q^{r-1}}) \]

	See Fulman \cite{fulcycle} for details.

\item Fulman \cite{fulthesis} develops probabilistic algorithms for growing
the random partitions $\lambda_{\phi}$ (after one chooses the size $n$
randomly according to a geometric distribution). It would be splendid if
these (or other) algorithms could be used to give a probabilistic proof of
Rogers-Ramanujan.

\item Ian Macdonald remarked to the author that the sum sides of the two
Rogers-Ramanujan identities

\[ \sum_{n=0}^{\infty} \frac{x^{n^2}}{(1-x)(1-x^2) \cdots (1-x^n)} =
\prod_{m=1}^{\infty} \frac{1}{(1-x^{5m-4})(1-x^{5m-1})} \]

\[ \sum_{n=0}^{\infty} \frac{x^{n(n+1)}}{(1-x)(1-x^2) \cdots (1-x^n)} =
\prod_{m=1}^{\infty} \frac{1}{(1-x^{5m-3})(1-x^{5m-2})} \]
 
	become (upon setting $x=\frac{1}{q}$)

\[ \sum_{n=0}^{\infty} \frac{1}{|GL(n,q)|} \]

\[ \sum_{n=0}^{\infty} \frac{1}{|AGL(n,q)|} \]

	where $AGL(n,q)$ is the affine finite general linear group. This
suggests that an analog of Theorem \ref{Rogers} should exist for the affine
finite general linear groups.

\item Is there some group theoretic reason why the right hand side of
Theorem \ref{Rogers} is almost a modular form? See Kac's beautiful article
on modular invariance \cite{Kac}.

\end{enumerate}

\section{Connection with Symmetric Function Theory} \label{SYM}

	Theorem \ref{Rogers} has a very clean statement in terms of
symmetric functions. For this, we recall the Hall-Littlewood polynomials
(page 208 of Macdonald \cite{Mac}). Let the permutation $w$ act on
variables $x_1,x_2,\cdots$ by sending $x_i$ to $x_{w(i)}$. There is also a
coordinate-wise action of $w$ on $\lambda=(\lambda_1, \cdots,\lambda_n)$
and $S^{\lambda}_n$ is defined as the subgroup of $S_n$ stabilizing
$\lambda$ in this action. Recall that $m_i(\lambda)$ is the number of parts
of $\lambda$ of size $i$. For a partition
$\lambda=(\lambda_1,\cdots,\lambda_n)$ of length $\leq n$, two definitions
of the Hall-Littlewood polynomials are:

\begin{eqnarray*}
P_{\lambda}(x_1,\cdots,x_n;t) & = & [\frac{1}{\prod_{i \geq 0}
\prod_{r=1}^{m_i(\lambda)} \frac{1-t^r}{1-t}}] \sum_{w \in S_n}
w(x_1^{\lambda_1} \cdots x_n^{\lambda_n} \prod_{i<j} \frac{x_i-tx_j}
{x_i-x_j})\\
& = & \sum_{w \in S_n/S_n^{\lambda}} w(x_1^{\lambda_1} \cdots x_n^{\lambda_n} \prod_{\lambda_i > \lambda_j} \frac{x_i-tx_j}{x_i-x_j})
\end{eqnarray*}

	At first glance it is not obvious that these are polynomials, but
the denominators cancel out after the symmetrization. The Hall-Littlewood
polynomials interpolate between the Schur functions ($t=0$) and the
monomial symmetric functions ($t=1$).

	The statement of Theorem \ref{symint} uses the standard notation
that $n(\lambda) = \sum_{i \geq 1} (i-1) \lambda_i$.

\begin{theorem} \label{symint} For $k \geq 2$,

\[ \sum_{\lambda: \lambda_1<k}
\frac{P_{\lambda}(\frac{1}{q},\frac{1}{q^2},\cdots;\frac{1}{q})}
{q^{n(\lambda)}} = \prod_{r=1 \atop r \neq 0, \pm k (mod \ 2k+1)}^{\infty}
(\frac{1}{1-\frac{1}{q^r}}) \]

\end{theorem}

\begin{proof}
	The result follows from Theorem \ref{Rogers} and Macdonald's
principal specialization formula (page 337 of Macdonald \cite{Mac}) which,
when applied to the Hall-Littlewood polynomials, states that

\begin{eqnarray*}
\frac{P_{\lambda}(\frac{1}{q},\frac{1}{q^2},\cdots;\frac{1}{q})}
{q^{n(\lambda)}} & = &  \frac{1}{q^{|\lambda|+2n(\lambda)} \prod_i
(\frac{1}{q})_{m_i(\lambda)}}\\
& = & \frac{1}{q^{\sum_i (\lambda_i')^2} \prod_i
(\frac{1}{q})_{m_i(\lambda)}}
\end{eqnarray*}

\end{proof}

{\bf Remarks}

\begin{enumerate}

\item Stembridge \cite{Stembridge} also found some connections between the
Rogers-Ramanujan identities and the Hall-Littlewood symmetric functions,
but the statement of Theorem \ref{symint} seems to be new.

\item The Schur functions are well known to give rise to the irreducible
polynomial representations of the general linear groups (see for instance Chapter
3 of Macdonald \cite{Mac}. Theorem
\ref{symint}, together with the paper of Feigen and Frenkel \cite{Feigen},
suggests that the Hall-Littlewood polynomials should also have a
representation theoretic interpretation.

\end{enumerate}
 
\section{Acknowledgements} This work is taken from the author's
Ph.D. thesis, done under the supervision of Persi Diaconis. His idea of
studying the random partitions $\lambda_{\phi}$ led to this work. Ed
Frenkel and Jim Lepowsky provided helpful pointers to the Rogers-Ramanujan
literature. This research was done under the generous 3-year support of the
National Defense Science and Engineering Graduate Fellowship (grant
no. DAAH04-93-G-0270) and the support of the Alfred P. Sloan Foundation
Dissertation Fellowship.

\end{document}